\documentclass[12pt]{amsart}

\usepackage{amsmath,amsfonts,amssymb}

\theoremstyle{plain}
\newtheorem{lemma}{Lemma}[section]
\newtheorem{theorem}[lemma]{Theorem}
\newtheorem{corollary}[lemma]{Corollary}

\theoremstyle{definition}

\numberwithin{equation}{section}
\begin{document}

\newcommand{\ZZ}{\mathbb{Z}}
\newcommand{\ZZd}{\mathbb{Z}^{d}}
\newcommand{\RR}{\mathbb{R}}
\newcommand{\RRd}{\mathbb{R}^{d}}
\newcommand{\PP}{\mathbb{P}}
\newcommand{\QQ}{\mathbb{Q}}
\newcommand{\EE}{\mathbb{E}}
\newcommand{\mA}{\mathcal{A}}
\newcommand{\mB}{\mathcal{B}}
\newcommand{\mC}{\mathcal{C}}
\newcommand{\mD}{\mathcal{D}}
\newcommand{\mE}{\mathcal{E}}
\newcommand{\mF}{\mathcal{F}}
\newcommand{\mG}{\mathcal{G}}
\newcommand{\mH}{\mathcal{H}}
\newcommand{\mJ}{\mathcal{J}}
\newcommand{\mL}{\mathcal{L}}
\newcommand{\mkN}{\mathfrak{N}}
\newcommand{\mR}{\mathcal{R}}
\newcommand{\mS}{\mathcal{S}}
\newcommand{\mT}{\mathcal{T}}
\newcommand{\mU}{\mathcal{U}}
\newcommand{\mW}{\mathcal{W}}
\newcommand{\bs}{\backslash}
\newcommand{\half}{\frac{1}{2}}
\newcommand{\fN}{\frac{1}{N}}
\newcommand{\tP}{\tilde{P}}
\newcommand{\hP}{\hat{P}}
\newcommand{\tf}{\tilde{f}}
\newcommand{\hu}{\hat{u}}
\newcommand{\hd}{\hat{\delta}}
\newcommand{\hf}{\hat{f}}

\title[First-Order Depinning Transitions]
{Ivy on the Ceiling:  First-Order Polymer Depinning Transitions with Quenched Disorder}
\author{Kenneth S. Alexander}
\address{Department of Mathematics KAP 108\\
University of Southern California\\
Los Angeles, CA  90089-2532 USA}
\email{alexandr@math.usc.edu}
\thanks{Research supported by NSF grant DMS-0405915.}

\date{\today}

\begin{abstract}
We consider a polymer, with monomer locations modeled by the
trajectory of an underlying Markov
chain, in the presence of a potential that
interacts with the polymer when it visits a particular site 0.  Disorder is
introduced by having the interaction vary from one monomer to
another, as a constant $u$ plus i.i.d.\ mean-0 randomness.  There is a
critical value
of $u$ above which the polymer is pinned, placing a positive fraction
(called the contact fraction) of its monomers
at 0 with high probability.  When the excursions of the underlying chain have a finite mean but no finite exponential moment, it is known \cite{AS06} that the depinning transition (more precisely, the contact fraction) in the corresponding annealed system is discontinuous.  One generally expects the presence of disorder to smooth transitions, and it is known \cite{GT06a} that when the excursion length distribution has power-law tails, the quenched system has a continuous transition even if the annealed system does not.  We show here that when the underlying chain is transient but the finite part of the excursion length distribution has exponential tails, then the depinning transition is discontinuous even in the quenched system, and the quenched and annealed critical points are strictly different.  By contrast, in the recurrent case, the depinning behavior depends on the subexponential prefactors on the exponential decay of the excursion length distribution, and when these prefactors decay with an appropriate power law, the quenched transition is continuous even though the annealed one is not. 
\end{abstract}
\maketitle

\section{Introduction} \label{S:intro}
It is a well-established principle in statistical mechanics that quenched disorder tends to smooth phase transitions.  One rigorous version of this principle was recently proved by Giacomin and Toninelli \cite {GT06a}, for polymer models of the following type.  The configuration of the polymer (in the absence of a potential) is described by the space-time trajectory of a Markov chain, which we take to be aperiodic; the location of the $i$th monomer is the state the chain is in, at time $i$.  A potential is then added at one site 0 of the state space (or in a one-dimensional subspace, from the perspective of space-time trajectories), attracting or repelling the polymer whenever site 0 is visited.  Quenched disorder can be incorporated by allowing the attraction of the potential to vary randomly from monomer to monomer, taking the form $u+V_i$ for the $i$th monomer, with $\{V_i\}$ i.i.d.\ with mean 0.  The corresponding Gibbs weight of a trajectory 
$x_{[0,N]} = \{ x_i: 0 \leq i \leq N \}$ is
\begin{equation} \label{E:GibbsWtQ}
    W_{\beta,u,P^X}(x_{[0,N]}) = \exp\left( \beta \sum_{i=1}^N (u + V_i)\delta_{\{x_i = 0 \}} \right)
          P^X(x_{[0,N]}),
\end{equation}
where $P^X$ denotes probability for the underlying Markov chain $\{X_i\}$, and $\beta$ denotes the inverse temperature.  
In \cite{GT06a} it is proved that for a large class of such models, the free energy rises at most quadratically in $\Delta = u-u_c^q$ as the pinning potential $u$ increases from its (quenched) critical point $u_c^q = u_c^q(\beta)$.  In other words, the quenched specific heat exponent (roughly speaking, the value $\alpha$ such that the free energy grows like $\Delta^{2-\alpha}$) is non-positive.
The fraction of monomers in contact with the pinning potential, called the \emph{contact fraction}, is the derivative in $u$ of the free energy, and hence rises at most linearly in $\Delta$; in particular, the contact fraction is continuous at $u_c^q$.  By contrast, if the (possibly infinite) excursion length for the underlying Markov chain lacks a finite exponential moment, but has a finite mean when conditioned to be finite, then in the annealed system (or equivalently, in the system with $V_i \equiv 0$) the contact fraction is discontinuous \cite{AS06}.
The key property assumed in \cite{GT06a} is that the excursion length distribution for the underlying Markov chain has power-law tails: there exists a finite $c$ such that for an excursion length $E$,
\begin{equation} \label{E:polytail}
  P^X(E=n) \geq n^{-c} \quad \text{ for all } n.
  \end{equation}
Suppose that in fact, for some $c \geq 1$ and slowly varying function $\varphi$,
\begin{equation} \label{E:tails}
      P^X(E = n) = n^{-c}\varphi(n).
\end{equation}
Then, up to slowly varying factors, the annealed free energy is proportional to $\Delta^{1/(c-1)}$ and the annealed contact fraction is proportional to $\Delta^{(2-c)/(c-1)}$ (\cite{Al06}, \cite{GT06a}).  (This time, $\Delta$ represents the increment from the annealed critical point.)  Thus the result in \cite{GT06a} is consistent with the Harris criterion, which in one form states that disorder is ``relevant,'' i.e. alters the specific heat exponent, when this exponent is positive ($c>3/2$), but is irrelevant when this exponent is negative ($c<3/2$.)  The Harris criterion is in general a nonrigorous principle from the physics literature, but a related result is established rigorously in \cite{CCFS89}.

Note that the conditions for a discontinuous transition in the annealed system are satisfied whenever \eqref{E:tails} holds with $c>2$.

The smoothing of the depinning transition due to disorder can be viewed heuristically as follows.  When the mean potential $u$ is near the critical point, due to random fluctuations there will be ``good'' stretches of the disorder in which the average of the $u+V_i$'s  is above $u_c$, and ``bad'' stretches where it is below.  The polymer can configure itself so that the good stretches are pinned, and the bad stretches are not.  As $u$ increases, the good stretches become scarcer, and the polymer gradually depins.
  
There are interesting cases in which \eqref{E:polytail} is not satisfied.  For example, one can model a force pulling the polymer away from the potential by giving the Markov chain a drift away from 0.  Let us consider in particular a biased simple random walk $\{X_i\}$ on the integers, with $P^X(X_{i+1} = x+1 \mid X_i = x) = p > 1/2, P^X(X_{i+1} = x-1 \mid X_i = x) = 1-p$, conditioned to stay nonnegative, with a random potential at state 0.  We denote the distribution of the unconditioned random walk here by $P^X_p$.  Without the potential, the walk is transient and the finite part of the excursion length distribution has an exponential tail.  In the annealed case, the transition is therefore first order \cite{AS06}.  We will show here that in this model and others like it, the transition remains first order in the quenched case.  

As described below (see \eqref{E:GibbsForce}), the model \eqref{E:GibbsWtQ} with $P^X = P_p^X$ is equivalent to one in which trajectories (of length $2N$) are confined to the upper half plane, $P^X$ is $P_{1/2}^X$ conditioned on the trajectory staying in the upper half plane, and there is an additional factor of $(p/(1-p))^{X_{2N}/2}$ in the Gibbs weight.  This additional factor reflects a force pulling the polymer away from the surface to which it is pinned (see \cite{Gi06}, \cite{GT06c}, \cite{IOW04}).

Besides changes in the order of the transition and in the specific heat exponent, it is of interest to know whether the critical point differs between the quenched and annealed systems.  The question is of interest in part because it is intertwined with questions of just how the polymer depins as the quenched critical point is approached, or, put differently, questions of what ``strategy'' the polymer uses to stay pinned when near the critical point--see \cite{CGG06}.  Simulations and nonrigorous methods \cite{NN01} have suggested shifts in the critical point in some cases where the excursion length has a power-law tail, but there are no rigorous proofs, and it was proved in \cite{Al06} that the critical points are the same in the case of \eqref{E:tails} with $c<3/2$, or $c=3/2$ with $\sum_n n^{-1} \varphi(n)^{-2} < \infty$.  We will show that in the above biased simple random walk case, and more generally whenever the underlying chain is transient and the finite part of the excursion length distribution has an exponential tail, the quenched critical point is strictly larger.  

The heuristic for discontinuous depinning in the quenched case is as given in the title: ivy on the ceiling.  Imagine an ivy vine along which there are a sequence of sites of varying stickiness, by way of which it can adhere to a ceiling despite the downward pull of gravity.  In winter the cold reduces the stickiness and an increasing number of sites let go of the ceiling.  At some point the number of attached sites becomes insufficient and the ivy lets go ``catastrophically.''  For a sufficiently long vine, the fraction of attached sites will not decrease to 0 as the stickiness decreases to its critical value (as would be the case in the absence of gravity); instead, there should be a certain strictly positive minimum contact fraction needed to counter gravity and keep the ivy attached.  Thus the transition should be first order.  A similar catastrophic-depinning heuristic (this time varying the force rather than the pinning potential) suggests why bandages often release suddenly and painfully when pulled from skin, though the physics of adhesives are of course more complex than the depinning phenomenon studied here.

\section{Results}
We use $P^V$ and $\langle \cdot \rangle^V$, respectively, to denote probability and expectation for the disorder $\{V_i\}$, and $\langle \cdot \rangle^Q$ and $E^Q(\cdot \mid \cdot)$, respectively, to denote expectation and conditional expectation for the chain $\{X_i\}$ under a measure $Q$.  As given, \eqref{E:GibbsWtQ} defines the \emph{quenched} version of the polymer model; the \emph{annealed} version is obtained by replacing the Gibbs weight \eqref{E:GibbsWtQ} with its $P^V$-expectation.  This is equivalent to taking $V_i \equiv 0$ and replace $u$ with $u + \beta^{-1} \log M_V(\beta)$, where $M_V$ is the moment generating function of $V_1$.  The corresponding annealed partition function and finite-volume Gibbs measure are denoted $Z_N(\beta,u)$ and $\mu_N^{\beta,u}$, respectively.  Let
\[
        L_N = \sum_{i=1}^N \delta_{\{X_i = 0 \}},
\]
where $\delta_A$ denotes the indicator of the event $A$, and let $E_i \leq \infty$ denote the length of the $i$th excursion from 0 for the chain $\{X_i\}$.  
The annealed free energy $f^a(\beta,u)$ is given by
\[
  \beta f^a(\beta,u) =  \lim_N \fN \log Z_N(\beta,u).
  \]
The \emph{contact fraction} for the annealed system is the unique value $C = C^a(\beta,u)$ for which 
\begin{equation} \label{E:contfrac}
     \lim_{N \to \infty} \mu^{\beta,u}_N\left( \frac{L_N}{N} \in
      (C-\epsilon,C+\epsilon) \right) = 1 \quad \text{for all } \epsilon > 0;
\end{equation}
the existence of such a $C$ is established in \cite{AS06}.  
Now $f^a$ is a convex
functions of $u$, and we have by standard methods that
\[
      C^a(\beta,u) = \frac{\partial}{\partial u} f^a(\beta,u)
\]
for all non-critical $u$.  The necessary differentiability of $f^a(\beta,\cdot)$ here follows from the fact that for supercritical $u$, $\beta f^a(\beta,u)$ is the solution $x>0$ of 
\[
  \sum_{1 \leq n < \infty} P^X(E_1 = n)e^{-xn} = \frac{ e^{-\beta u} }{ M_V(\beta) };
  \]
see Appendix A of \cite{GT06a}.

Let us call a Markov chain \emph{nontrivially transient} if $0<P^X(E_1 < \infty)<1$.  We say that the (annealed) polymer is \emph{pinned} at $(\beta,u)$ if the contact fraction is positive.  The annealed critical point is 
\[
  u_c^a = u_c^a(\beta) = \inf\{u \in \RR: C^a(\beta,u)>0\},
  \]
which from \cite{AS06} is in $[-\infty,\infty)$ provided the chain is not trivially transient.
  
For the quenched system, the partition function and finite-volume Gibbs measure are denoted $Z_N^{ \{V_i\} }(\beta,u)$ and $\mu_N^{\beta,u,\{V_i\}}$ respectively, and we have definitions analogous to the annealed case for the free energy $f^q(\beta,u)$ and the quenched critical point $u_c^q(\beta)$.  In \cite{AS06} it was established that self-averaging holds in the sense that the quenched free energy exists and is nonrandom, provided we exclude a $P^V$-null set.  Differentiability of $f^q(\beta,\cdot)$ is proved in \cite{GT06b} when the underlying Markov chain has power-law tails on the excursion length distribution, but is not known in general, forcing us to define $C^{q,-}(\beta,u)$ and $C^{q,+}(\beta,u)$ to be the left and right derivatives respectively of the convex function $f^q(\beta,\cdot)$ at $u$, and then define 
\[
  \mD(\beta) = \{u \in \RR: u_c^q(\beta): C^{q,+}(\beta,u) = C^{q,-}(\beta,u)\},
  \]
a set for which the complement is at most countable.  For $u \in \mD(\beta)$ we denote the common value $\frac{\partial}{\partial u} f^q(\beta,u)$ by $C^q(\beta,u)$.  
From convexity and monotonicity of $f^q$ in $u$ we have for fixed $\beta$ that 
\[
  \left\langle \frac{L_N}{N} \right\rangle^{\beta,u,\{V_i\}}_{[0,N]} = \frac{1}{\beta} \frac{\partial}{\partial u} \left(
    \fN \log Z^{\{V_i\}}_{[0,N]}(\beta,u) \right) \to \frac{\partial}{\partial u} 
    f^q(\beta,u) \text{ for all } u \in \mD(\beta)
  \]
and 
\[
  C^{q,-}(\beta,u) \leq \liminf_N \left\langle \frac{L_N}{N} \right\rangle^{\beta,u,\{V_i\}}_{[0,N]}
    \leq \limsup_N  \left\langle \frac{L_N}{N} \right\rangle^{\beta,u,\{V_i\}}_{[0,N]} \leq C^{q,+}(\beta,u) 
    \text{ for all } u \notin \mD(\beta), 
    \]
both $P^V$-a.s.
Let 
\[
  b_E = b_E(P^X) = \sup\{b \geq 0: P^X(E_1 = n) =O(e^{-bn}) \text{ as } n \to \infty\}.
  \]
When $b_E(P^X)>0$ we say that the chain, or $P^X$, has \emph{exponential excursion tails}.  Note this does not rule out the possibility that the chain is transient.    When $b_E < \infty$ we define the prefactors $\gamma_n$ by
\[
  P^X(E_1 = n) = \gamma_n e^{-b_E n}.
  \]

It is shown in \cite{AS06} that for a nontrivially transient chain with exponential excursion tails, $u_c^a(\beta)$ and $u_c^q(\beta)$ are finite for all $\beta>0$.

\begin{theorem} \label{T:main}
Suppose that $\{V_i, i \geq 1\}$ are i.i.d.\ Gaussian random variables, and the underlying Markov chain is nontrivially transient with exponential excursion tails.  Then
\begin{itemize}
\item[(i)] the quenched and annealed transitions are both discontinuous, in that $C^q(\beta,u)$ and $C^a(\beta,u)$ are discontinuous in $u$ at $u_c^q(\beta)$ and $u_c^a(\beta)$, respectively;
\item[(ii)] the two critical points are strictly different:  $u_c^a(\beta) < u_c^q(\beta)$ for all $\beta>0$.
\end{itemize}
\end{theorem}

The assumption in Theorem \ref{T:main} that the disorder is Gaussian is used only to allow the citation of a result from \cite{Al06} (see \eqref{E:compcost} below) which is only proved in \cite{Al06} for the Gaussian case.  In the case of general disorder having a finite exponential moment, the term $\half \beta^2 \delta^2$ in \eqref{E:compcost} becomes only an approximation valid for small $\delta$, but this should affect only the technical details, so we expect that Theorem \ref{T:main} is valid for arbitrary disorder distributions having a finite exponential moment.

Consider a random walk trajectory $\{x_{[0,2N]}\} \in \ZZ^{2N+1}$, and let $p \neq 1/2 \in (0,1)$.  We have
\begin{equation} \label{E:PXpform}
  P_p^X(x_{[0,2N]}) = p^{N+x_{2N}/2}(1-p)^{N-x_{2N}/2} 
    = (4p(1-p))^N P_{1/2}^X(x_{[0,2N]}) \left( \frac{p}{1-p} \right)^{x_{2N}/2},
    \end{equation}
so the measure $P_p^X$ (which has exponential excursion tails) is equivalent to the measure $P_{1/2}^X$ (which satisfies \eqref{E:tails}) weighted by $(p/(1-p))^{x_{2N}/2}$.  The same weighting applied to \eqref{E:GibbsWtQ} yields Gibbs weights
\begin{equation} \label{E:GibbsForce}
        \exp\left( \beta \sum_{i=1}^{2N} (u + V_i)\delta_{\{x_i = 0 \}} \right)
          P^X_{1/2}(x_{[0,2N]}) \left( \frac{p}{1-p} \right)^{x_{2N}/2}.
\end{equation}
We can further alter this by restricting to space-time trajectories in the upper half plane, replacing $P^X_{1/2}$ with $P^X_{1/2}( \{X_i\} \in \cdot \mid X_i \geq 0$ for all $i \leq 2N)$.  With respect to pinning behavior, this is equivalent to allowing trajectories in the full plane but multiplying the Gibbs weight of a trajectory with $k$ returns to 0 by  
$2^{-k}$, and this multiplication is further equivalent to replacing $u$ with $u - \beta^{-1} \log 2$.
The model with weights given by \eqref{E:GibbsForce} with trajectories restricted to the upper half plane is an instance of a model considered in \cite{Gi06}, \cite{GT06c}, \cite{IOW04} and elsewhere for a polymer depinned by a force pulling on the free end.  Thus for that model we immediately have the following.
\begin{corollary} \label{T:force}
Suppose that $\{V_i, i \geq 1\}$ are i.i.d.\ Gaussian random variables and $p > 1/2$.  Then the model with weights given by \eqref{E:GibbsForce} with trajectories restricted to the upper half plane satisfies (i) and (ii) of Theorem \ref{T:main}.
\end{corollary}

The case of a recurrent underlying Markov chain with exponential excursion tails, in the cases we can deal with, is generally an elementary application of existing results, but is useful as a contrast to the transient case.
Our first result for the recurrent case includes correction of an error in (\cite {AS06}, Theorem 2.1), where it was stated that in the case of a recurrent underlying chain with exponential excursion tails, the transition in the annealed system was necessarily continuous.  It is stated in the proof of that theorem that the function $g$, defined there and in the proof of Theorem \ref{T:main} below, is strictly convex on $[0,m_E^{-1}]$, but this is not necessarily true.  The corrected result below is written for the deterministic system ($V_i \equiv 0$), but is valid for the annealed system with disorder $V_i$ having a finite exponential moment, for $\beta$ for which $M_V(\beta) < \infty$, since the annealed system at $(\beta,u)$ is the same as the deterministic system at $(\beta,u + \beta^{-1}\log M_V(\beta))$.  For the deterministic system we omit the superscript $q$ or $a$ on the free energy, contact fraction, etc.

Let
\[
  M_E^f(x) = E^{P^X}\left( e^{xE_1} \mid E_1 < \infty \right), 
    \qquad b_E' = \lim_{x \nearrow b_E} (\log M_E^f)'(x),
  \]
or in terms of the prefactors,
\[
  M_E^f(b_E) = \sum_n \gamma_n, \qquad b_E' = \sum_n n\gamma_n \quad \text{if } b_E < \infty.
  \]
When $M_E^f(b_E) < \infty$ and $P^X$ is recurrent, we can define a measure $\tP^X$ by
\[
  \tP^X(E_1 = n) = \frac{\gamma_n}{M_E^f(b_E)}.
  \]
We complete the definition of $\tP$ by specifying that the distribution given the excursion lengths be the same as under $P^X$. 
We refer to the system with $P^X$ replaced by $\tP^X$ as \emph{loosened} since $\tP^X$ does not have exponential excursion tails, and we write the corresponding free energy in the deterministic case as $\tf(\beta,u)$.  When distinguishing it from the loosened system we will refer to the system under $P^X$ as \emph{original}.
Observe that when the underlying chain is biased random walk on $\ZZ$ (measure $P_p^X$), by \eqref{E:PXpform} we have $b_E = -\half \log (4p(1-p))$ and the loosened measure is  $\tP = P_{1/2}^X$.

\begin{theorem} \label{T:corrected}
Suppose that $V_i \equiv 0$ and the underlying Markov chain is recurrent with exponential excursion tails.  Then
\begin{itemize}
\item[(i)] if $M_E(b_E) = \infty$ then there is no transition (that is, $u_c(\beta) = -\infty$ for all $\beta>0$);
\item[(ii)] if $M_E(b_E) < \infty$ and $b_E' = \infty$ then there is a continuous transition at $u_c(\beta) = -\beta^{-1} \log M_E(b_E)$ for all $\beta>0$ (that is, $C(\beta,u)$ is continuous in $u$ at $u_c(\beta)$);
\item[(iii)] if $M_E(b_E) < \infty$ and $b_E' < \infty$ then there is a discontinuous transition at $u_c(\beta) = -\beta^{-1} \log M_E(b_E)$ for all $\beta>0$.
\end{itemize}
Further, if $M_E(b_E) < \infty$ then
\begin{equation} \label{E:freeenrel}
  \beta f(\beta,u) = \beta \tf(\beta,u + \beta^{-1}\log M_E(b_E)) - b_E \quad \text{for all } \beta, u.
\end{equation}
\end{theorem}

Equation \eqref{E:freeenrel} says that the graph of $f(\beta, \cdot)$ is just a translate of the graph of $\tf(\beta,\cdot)$, so the nature of the transition is the same in the original and loosened systems.

Giacomin and Toninelli \cite{GT06a} showed that if 
\[
  P^X(E_1 = n) \geq n^{-c} \quad \text{for all $n \geq 1$, for some } c \geq 1,
  \]
and $V_1$ is either bounded or has a density $\varphi$ with respect to Lebesgue measure satisfying
\begin{equation} \label{E:GTcond2}
  \int_{\RR} \varphi(x+y) \log\left( \frac{ \varphi(x+y) }{ \varphi(y) } \right)\ dy \leq Rx^2 
    \quad \text{ for all sufficiently small } x
    \end{equation}
for some $R>0$, then the transition is continuous in the quenched system.  We can apply this to the case of exponential excursion tails as follows, showing that with exponential excursion tails, the recurrent case is quite different from the transient case.

\begin{theorem} \label{T:recurquen}
Suppose that $\{V_i, i \geq 1\}$ are i.i.d.\ and the underlying Markov chain is recurrent with exponential excursion tails, satisfying $M_E(b_E) < \infty$ and
\begin{equation} \label{E:regdecay}
  \lim_n \frac{1}{n} \log P^X(E_1 > n) = -b_E.
  \end{equation}
Then
\begin{equation} \label{E:freeenrel2}
  \beta f^q(\beta,u) = \beta \tf^q(\beta,u + \beta^{-1}\log M_E(b_E)) - b_E \quad \text{for all } \beta, u,
\end{equation}
and therefore the transition is continuous in the quenched system if and only if it is continuous in the loosened quenched system.  In particular if the prefactors satisfy
\[
  \gamma_n \geq n^{-c} \quad \text{for all $n \geq 1$, for some } c \geq 1,
  \]
and $V_1$ is either bounded or has a density $\varphi$ with respect to Lebesgue measure satisfying \eqref{E:GTcond2}, then the transition is continuous in the quenched system.
\end{theorem}

In the recurrent case in which the prefactors $\gamma_n$ decay faster than any power of $n$, the annealed system has a discontinuous transition but we do not know whether the transition is continuous in the quenched system.  Theorem \ref{T:recurquen} only turns the question into a similar one about the loosened system.

\section{Proofs}
We will give two proofs of Theorem \ref{T:main} which provide quite different intuition, as each may be useful in other contexts.  The second proof will follow the proofs of Theorems \ref{T:corrected} and \ref{T:recurquen}; here is the first.

\begin{proof}[Proof of Theorem \ref{T:main}]
Discontinuity of the transition in the annealed system was established in \cite{AS06}, so we consider the quenched system.  We may assume the underlying chain is aperiodic, and we need only consider $u \in \mD(\beta)$.
Write $u$ as $u = u_c^a + \Delta$ with $\Delta \geq 0$.  From \cite{AS06} we have $u_c^q(\beta) < u_c^a(\beta) + \beta^{-1} \log M_V(\beta)$ (the latter being the critical point for the deterministic system) so it suffices to consider $0 < \Delta < \beta^{-1} \log M_V(\beta)$.

Let
\[
     I_E(t) = -\lim_{\epsilon \searrow 0}\
       \lim_n \frac{1}{n} \log P^X \left( \frac{1}{n}\sum_{i=1}^n E_i
       \in (t-\epsilon,t+\epsilon) \right), \quad t>0,
\]
denote the large-deviations rate function for $E_1$, let $r = -\log P^X(E_1 < \infty)$, let
\begin{align}
  I_E^f(t) &= I_E(t) - r \notag \\
  &= -\lim_{\epsilon \searrow 0}\
       \lim_n \frac{1}{n} \log P^X \left( \frac{1}{n}\sum_{i=1}^n E_i
       \in (t-\epsilon,t+\epsilon)\ \bigg|\ E_1 < \infty,..,E_n < \infty \right), \notag \\
   &= \sup_x (tx - \log M_E^f(x)), \notag
\end{align}
and let
\[
  M_E(x) = \left\langle e^{xE_1} \right\rangle^{P^X} \ \  (x \neq 0), \quad M_E(0) = P^X(E_1 < \infty),
    \]
so $M_E$ is left-continuous at 0.
Let $m_E = E^X(E_1 \mid E_1 < \infty)$, and let
\[
  a = \min\{n: P^X(E_1 = n)>0\}, \quad A = \sup\{n<\infty: P^X(E_1 = n) > 0\}.
  \]
Since (i) is trivial for $A<\infty$, we assume for now that $A=\infty$.  This means that if $b_E = \infty$ then $b_E' = \infty$.
Note that by convexity $b_E' > (\log M_E^f)'(0) = m_E$.  Let 
\[
   J_E(t) = \sup_x \left( tx - \log M_E(x) \right), \quad t>0,
\]
and let
\[
  g(x) = \begin{cases} xJ_E(x^{-1}), &x \in (0,1]\\ 0, &x=0, \end{cases}
  \]
\[
  \hat{g}^f(x) = \begin{cases} xI_E^f(x^{-1}), &x \in (0,1],\\
    b_E, &x=0, \end{cases}
    \]
\[
  \hat{g}(x) = \hat{g}^f(x) + rx.
  \]
Of course in the recurrent case we would have $J_E = I_E$ but they differ here:  for all $t > m_E$ we have $-\log P^X(E_1 < \infty) = -\log M_E(0) = J_E(t) < I_E(t)$.  Basic properties include the following: 
\begin{itemize}
\item[(i)] $\hat{g}^f$ is infinite outside $[0,a^{-1}]$, has minimum $\hat{g}^f(m_E^{-1}) = 0$, is convex on $(0,a^{-1}]$, is affine on $(0,(b_E')^{-1}]$ and is strictly convex on $[(b_E')^{-1},a^{-1}]$, 
\item[(ii)] $\hat{g}$ is minimized at some $x^* \in [0,m_E^{-1})$, 
\item[(iii)] $g, \hat{g}, \hat{g}^f$ are continuous (as extended-real-valued functions) at 0,
\item[(iv)] $\lim_{x \searrow 0}(\hat{g}^f)'(x) = -\log M_E^f(b_E) \in [-\infty,\infty)$,
\item[(v)] $g$ is the convex minorant of the function 
\[
  g_0(x) = \begin{cases} \hat{g}(x) &\text{if } x>0\\ 0 &\text{if } x=0, \end{cases}
  \]
satisfying $g(x) = rx$ for $0 \leq x \leq m_E^{-1}$, $g(x) = \hat{g}(x)$ for $x \geq m_E^{-1}$.  
\end{itemize}
Here (i) follows from the fact that the infimum in the variational formula for $I_E^f(t)$ is achieved at $x=b_E$ if and only if $M_E(b_E) < \infty$ and $t \geq b_E'$, and then $I_E^f$ is affine for such $t$; otherwise this infimum is achieved in $(-\infty,b_E)$ where $\log M_E$ is strictly convex and analytic.
The other nontrivial facts among (i)--(v) were established in (\cite{AS06}, proof of Lemma 2.2).  In (v) the value $g_0(0) = 0$ corresponds to the fact that transience allows the event of ``no returns to 0'' to occur at a cost which is constant, so in particular does not decay exponentially in $N$.
From \cite{AS06}, first we have that for all $0 \leq \delta < \eta \leq 1$, 
\begin{equation} \label{E:LNdev}
      \lim_N  \fN \log P^X(\delta N < L_N < \eta N) =
      -\inf_{\delta < x < \eta} g(x),
\end{equation}
and second, noting that $r+\beta \Delta = \beta u + \log M_V(\beta)$, we have the variational principle
\begin{equation} \label{E:anfreeen}
     \beta f^a(\beta,u) = 
       \sup_{\delta \in [0,1]} ((r+\beta \Delta) \delta - g(\delta)).
\end{equation}
The contact fraction $C^a(\beta,u)$ is the value of $\delta$ which achieves this supremum.  This value is unique for $\Delta > 0$ by strict convexity of $g$--see (i) and (v) above.   In particular, 
\[
  \beta f^a(\beta,u) = 
       \max\left( \sup_{\delta \in [0,1]} \left((r+\beta \Delta) \delta - \hat{g}(\delta) \right), 0 \right).
       \]

We claim that if we consider only trajectories with $X_N = 0$, then we replace $g$ with $\hat{g}$ in \eqref{E:LNdev}, that is,
\begin{equation} \label{E:hatglimit}
      \lim_N  \fN \log P^X(\delta N < L_N < \eta N, X_N = 0) =
      -\inf_{\delta < x < \eta} \hat{g}(x).
\end{equation}
Let $\hat{g}_{\min}$ denote the infimum on the right side of \eqref{E:hatglimit}, and let $\tau_n = \sum_{i=1}^n E_i$ be the time of the $n$th return to 0, with $\tau_n = \infty$ if there is no such return.  Then
\begin{align} \label{E:upper1}
  P^X(\delta N < L_N < \eta N, X_N = 0) &= \sum_{\delta N < k < \eta N} P^X(\tau_k = N) \\
  &\leq \sum_{\delta N < k < \eta N} e^{-kI_E(N/k)} \notag \\
  &\leq Ne^{-\hat{g}_{\min}N}. \notag
  \end{align}
For a lower bound we need several cases.

\emph{Case 1}. $\delta < \eta \leq x^*$, for $x^*$ from (ii), with  $\eta > (b_E')^{-1}$.  Let $\epsilon>0$ satisfy $\eta^{-1} + 2\epsilon < \min(\delta^{-1},b_E')$, and define $\alpha = \alpha(\epsilon)$ by 
\begin{equation} \label{E:alphaprop}
  (\log M_E^f)'(\alpha) = \eta^{-1} + \epsilon.  
  \end{equation}
Since $(\log M_E^f)'(0) = m_E < (x^*)^{-1} \leq \eta^{-1}$, we must have $\alpha>0$.
Define the ``tilted'' measure $Q$ by
\[
  Q(E_1 = n) = \frac{ e^{r+\alpha n} P^X(E_1 = n) }{ M_E^f(\alpha) },
  \quad 1 \leq n < \infty.
  \]
We have $I_E(\eta^{-1} + \epsilon) = (\eta^{-1} + \epsilon) \alpha - \log M_E^f(\alpha) + r$, so
\begin{align} \label{E:lower1}
  P^X&(\delta N < L_N < \eta N, X_N = 0) \\
  &= \sum_{\delta N < k < \eta N} P^X(\tau_k = N) \notag \\
  &= \sum_{\delta N < k < \eta N} \exp\left( -\left( \frac{N}{k}\alpha - 
    \log M_E^f(\alpha) + r \right)k\right) Q(\tau_k = N) \notag \\
  &\geq \sum_{(\eta^{-1} + 2\epsilon)^{-1} N < k < \eta N} 
    \exp\left( -\left( \left( \eta^{-1} + 2\epsilon \right)\alpha - 
    \log M_E^f(\alpha) + r \right)k\right) Q(\tau_k = N) \notag \\
  &\geq \exp\left( - \eta \left( I_E(\eta^{-1} + \epsilon) -\epsilon\alpha \right) N \right)
    Q\big( \tau_k = N \text{ for some } (\eta^{-1} + 2\epsilon)^{-1} N < k < \eta N \big). \notag
  \end{align}
Since $\langle E_1 \rangle^Q = \eta^{-1} + \epsilon$, we can apply the renewal theorem (as in Appendix A of \cite{GT06a}) to obtain
\[
  Q\big( \tau_k = N \text{ for some } (\eta^{-1} + 2\epsilon)^{-1} N < k < \eta N \big) 
    \to \frac{1}{\eta^{-1} + \epsilon} \quad \text{as } N \to \infty.
    \]
Since $\epsilon$ is arbitrary, this and \eqref{E:upper1}, \eqref{E:lower1} show that 
\begin{equation} \label{E:ghatlimit}
      \lim_N  \fN \log P^X(\delta N < L_N < \eta N, X_N = 0) = -\hat{g}(\eta) = 
      -\hat{g}_{\min}.
\end{equation}

\emph{Case 2}. $\delta < \eta \leq x^*$, with  $\eta \leq (b_E')^{-1}$.   This means $b_E' < \infty$ so $b_E < \infty, \log M_E(b_E) < \infty$ and there is no $\alpha$ as in \eqref{E:alphaprop}, so instead we fix $0<\epsilon<(1 - \delta/\eta)/3$, take $\rho = (1 - 3\epsilon) b_E'$ and define $\alpha = \alpha(\epsilon)$ by $(\log M_E)'(\alpha) = (1 - 2\epsilon)^{-1}\rho$.  We have $I_E(\eta^{-1}) = \eta^{-1}b_E - \log M_E^f(b_E) + r$.  Since $b_E < \infty$ there exists a sequence $q_N \to \infty$ with
\[
  P^X(E_1 = q_N) = e^{-(b_E q_N + o(q_N))}.
  \]
Let $n_N = \max\{j \leq \eta\rho N:  N-j \text{ is a multiple of } q_N\}$ and $J_N = (N-n_N)/q_N$.
Defining the tilted measure $Q$ as in Case 1a, as in \eqref{E:lower1} we obtain provided $N$ is large that
\begin{align} \label{E:lower2}
  P^X&(\delta N < L_N < \eta N, X_N = 0) \\
  &\geq P^X(\text{for some } \delta N < k < (1-\epsilon)\eta N, \tau_k = n_N \text{ and } 
    E_{k+1} =..= E_{k+J_N} = q_N) \notag \\
  &\geq P^X(\tau_k = n_N \text{ for some } (1-3\epsilon)\rho^{-1}n_N < k < (1 - \epsilon)\eta N) 
    P^X(E_1 = q_N)^{J_N} \notag \\
  &\geq \exp\left( -\alpha n_N + \left(
    \log M_E^f(\alpha) - r \right) (1 - \epsilon) \eta N - b_E(N-n_N) - \epsilon N \right) \notag \\
  &\qquad \cdot Q\big( \tau_k = n_N \text{ for some } 
    (1-3\epsilon)\rho^{-1}n_N < k < (1-\epsilon)\rho^{-1} n_N \big). \notag
  \end{align}
Since $\langle E_1 \rangle^Q = (1 - 2\epsilon)^{-1}\rho$, the renewal theorem applies as above to the last probability, and since $\alpha \to b_E$ as $\epsilon \to 0$ we have from \eqref{E:upper1} and \eqref{E:lower2} that \eqref{E:ghatlimit} holds.

\emph{Case 3}.  $x^* \leq \delta < \eta$.  We may assume $\delta < a^{-1}$.  Let $\epsilon>0$ satisfy $\delta^{-1} - 2\epsilon > \max(\eta^{-1},a)$ and this time define $\alpha = \alpha(\epsilon)$ by $(\log M_E^f)'(\alpha) = \delta^{-1} - \epsilon$, then define $Q$ as in Case 1.  We have $I_E(\delta^{-1} - \epsilon) = (\delta^{-1} - \epsilon)\alpha - \log M_E^f(\alpha) + r$, so as in \eqref{E:lower1},
\begin{align} \label{E:lower3}
  P^X&(\delta N < L_N < \eta N, X_N = 0) \\
  &= \sum_{\delta N < k < \eta N} \exp\left( -\left( \frac{N}{k}\alpha - 
    \log M_E^f(\alpha) + r \right)k\right) Q(\tau_k = N) \notag \\
  &\geq \sum_{\delta N < k < (\delta^{-1}  - 2\epsilon)^{-1} N} 
    \exp\left( -\left( \left( \delta^{-1} - \epsilon \right)\alpha - 
    \log M_E^f(\alpha) + r + \epsilon|\alpha| \right)k\right) Q(T_k = N) \notag \\
  &\geq \exp\left( \left( -\delta (1 - 2\epsilon\delta)^{-1}
    I_E(\delta^{-1}  - \epsilon) - \epsilon |\alpha| \right) N \right) \notag \\
  &\qquad \cdot Q\big( \tau_k = N \text{ for some } \delta N < k < (\delta^{-1}  - 2\epsilon)^{-1} N \big).   
    \notag
  \end{align}
This time \eqref{E:upper1} is valid with $\delta$ in place of $\eta$ on the right side, so once again, since $\epsilon$ is arbitrary, \eqref{E:upper1}, the renewal theorem and \eqref{E:lower3} give
\begin{equation} \label{E:ghatlimit3}
      \lim_N  \fN \log P^X(\delta N < L_N < \eta N, X_N = 0) = -\hat{g}(\delta) = 
      -\hat{g}_{\min}.
\end{equation}

\emph{Case 4}. $\delta < x^* < \eta$.  We make use of Case 3:  for $\epsilon>0$ sufficiently small we have
\begin{align} \label{E:XN0cost}
  \liminf_N &\fN \log P^X(\delta N < L_N < \eta N, X_N = 0) \\
  &\geq \liminf_N \fN \log P^X(x^* N < L_N < \eta N, X_N = 0) \notag \\
  &= -\inf_{x^* < x < \eta} \hat{g}(x) \notag \\
  &=- \hat{g}(x^*) \notag \\
  &= -\hat{g}_{\min}. \notag
  \end{align}
Together with \eqref{E:upper1} this proves \eqref{E:hatglimit}, which is now proved in all cases.

Let $T_N$ denote the time of the last return to 0 in $[0,N]$.  
For $\Xi$ a set of trajectories of the chain, let $Z_k^{\{V_i\}}(\beta,u,\Xi)$ denote the contribution to the quenched partition function from trajectories in $\Xi$.  Then for $k \leq N$,
\begin{equation} \label{E:Zsplit}
  Z_N^{\{V_i\}}(\beta,u,\{ T_N = k \}) \leq Z_k^{\{V_i\}}(\beta,u,\{ X_k = 0\} )P^X(E_1 > N-k),
  \end{equation}
while from \eqref{E:hatglimit} and minor adaptations of the proof of (\cite{Al06}, Theorem 1.3) we have 
\begin{equation} \label{E:compcost}
  \limsup_k \frac{1}{k} \log Z_k^{\{V_i\}}(\beta,u,\{ X_k = 0\} ) \leq 
    \sup_{\delta \in [0,1]} \left( (r + \beta \Delta) \delta - \hat{g}(\delta) - \half \beta^2 \delta^2 \right).
  \end{equation}
Here we use the fact that $r + \beta\Delta = \beta u + \log M_V(\beta)$.  Considering $\lambda = k/N$ we see from \eqref{E:Zsplit} and \eqref{E:compcost} that
\begin{align} \label{E:qfreeineq}
  \beta f^q(\beta,u) &\leq
       \sup_{\lambda \in [0,1]} \left(
       \lambda \sup_{\delta \in [0,1]} \left( (r + \beta \Delta) \delta - \hat{g}(\delta) 
       - \half \beta^2 \delta^2 \right) \right) \notag \\
  &= \max\left( \sup_{\delta \in [0,1]} \left( (r + \beta \Delta) \delta - \hat{g}(\delta) 
    - \half \beta^2 \delta^2 \right), 0 \right).
  \end{align}
In fact, if $\theta > C^q(\beta,u)$ then we may replace $\sup_{\delta \in [0,1]}$ with $\sup_{\delta \in [0,\theta)}$ in \eqref{E:qfreeineq}.
Let
\[
  h(\delta) = \hat{g}(\delta) + \half \beta^2 \delta^2.
  \]
Then $h(0) = b_E \in (0,\infty]$, so there exists $0 < \delta_0 < \min(b_E,1)/4(\log M_V(\beta) + r)$ such that $\delta < \delta_0$ implies $h(\delta) > \min(b_E,1)/2$, and hence also implies
\[
  (r + \beta \Delta) \delta < \delta_0 (r + \log M_V(\beta)) < \frac{\min(b_E,1)}{4} < \frac{h(\delta)}{2}.
  \]
Thus in this case, $\sup_{\delta \in [0,\delta_0)} \left( (r + \beta \Delta) \delta - h(\delta) \right) < 0$.  From \cite{AS06}, we have $u>u_c^q(\beta)$ if and only if $\beta f^q(\beta,u) > 0$.  It follows that if $u>u_c^q(\beta)$ , then 
\[
  C^q(\beta,u) \geq \begin{cases} A^{-1}, &\text{if } A<\infty,\\ \delta_0, &\text{if } A=\infty, \end{cases}
  \]
proving that the transition in the quenched system is discontinuous.

Turning to the proof of (ii), we may assume $\delta_0 \leq m_E^{-1}$.  It is easily seen from \eqref{E:anfreeen} and (v) that $C^a(\beta,u) \geq m_E^{-1}$ for all $u > u_c^a(\beta)$.  Therefore we can restrict the sups in \eqref{E:anfreeen} and \eqref{E:qfreeineq} to $\delta \geq \delta_0$ for all $u$.  It follows that
\[
  \beta f^q(\beta,u) \leq \beta f^a(\beta,u) - \half \beta^2 \delta_0^2 \quad \text{for all } u \geq 
  u_c^q(\beta).
  \]
Since $f^a(\beta,u_c^a(\beta)) = 0$ and $\frac{\partial}{\partial u} f^a(\beta,u) = C^a(\beta,u) \leq 1$, it follows from this and \eqref{E:qfreeineq} that 
\[
  u_c^q(\beta) \geq u_c^a(\beta) + \half \beta \delta_0^2,
  \]
so the two critical points are distinct.
\end{proof}

\begin{proof}[Proof of Theorem \ref{T:corrected}]
Statement (i) is proved in \cite{AS06}, so we consider the case of $M_E(b_E) < \infty$.  We use tildes to denote quantities associated to the loosened system.  Then
\[
  \tilde{M}_E(x) = \frac{ M_E(x-b_E)}{M_E(b_E)},
  \]
so
\[
  \tilde{J}_E(t) = J_E(t) + tb_E - \log M_E(b_E),
  \]
and then
\[
  \tilde{g}(\delta) = g(\delta) +b_E - \delta \log M_E(b_E).
  \]
Now \eqref{E:freeenrel} is immediate from the variational principle \eqref{E:anfreeen}, so the transition in the original system is continuous if and only if the transition in the loosened system is continuous.  For the loosened system (ii) and (iii) are proved in \cite{AS06}, so we conclude they are also valid for the original system.
\end{proof}

\begin{proof}[Proof of Theorem \ref{T:recurquen}]
Let $\tilde{u} = u + \beta^{-1}\log M_E(b_E)$.
We have for a trajectory $x_{[0,N]}$ with $T_N = n$ for some $n \leq N$ that
\begin{equation} \label{E:weightratio}
  \frac{ W_{\beta,\tilde{u},\tP^X}(x_{[0,N]}) }{ W_{\beta,u,P^X}(x_{[0,N]}) } = 
    e^{b_E n} \frac{ \tP^X(E_1 > N-n) }{ P^X(E_1 > N-n) },
    \end{equation}
and \eqref{E:freeenrel2} follows easily from this and \eqref{E:regdecay}.
The rest of the theorem is immediate from \cite{GT06a}.
\end{proof}

\begin{proof}[Alternate proof of Theorem \ref{T:main}]
The loosened system need not exists in general, as we have not assumed $M_E(b_E) < \infty$, but we can construct a ``partially loosened'' system by defining a measure $\hat{P}$ as follows.  Let $0 < b < b_E$ and let $\hP$ have excursion length distribution
\[
  \hP^X(E_1 = n) = \frac{ P^X(E_1 = n)e^{bn}P^X(E_1 < \infty) }{ M_E^f(b) } \quad \text{for all } 1 \leq n < \infty,
  \]
\[ 
    \hP^X(E_1 = \infty) = P^X(E_1 = \infty).
    \]
We denote the quenched free energy, corresponding to the Gibbs weights $W_{\beta,u,\hP}(\cdot)$ from \eqref{E:GibbsWtQ}, by $\hf^q(\beta,u)$, the contact fraction by $\hat{C}^q(\beta,u)$, and the quenched critical point by $\hu_c^q(\beta)$.  In both quenched and annealed systems, and in both original and loosened, we have $u$ greater than the critical point if and only if the free energy is positive; otherwise the free energy is 0.
Now for $\hu = u + \beta^{-1}\log M_E^f(b)$ we have analogously to \eqref{E:weightratio} for a trajectory $x_{[0,N]}$ with $T_N = n$ that
\begin{equation} \label{E:weightratio2}
  \frac{ W_{\beta,\hu,\hP^X}(x_{[0,N]}) }{ W_{\beta,u,P^X}(x_{[0,N]}) } = 
    e^{bn} \frac{ \hP^X(E_1 > N-n) }{ P^X(E_1 > N-n) },
    \end{equation}
with the last ratio being bounded away from 0.
Now by \eqref{E:Zsplit}, for $\lambda < 1$,
\[
  \limsup_N \fN \log Z_N^{\{V_i\}}(\beta,u,\{ T_N \leq \lambda N \}) \leq \lambda \beta f^q(\beta,u),
  \]
so when $u > u_c^q(\beta)$, i.e. when $f^q$ is positive, we have $\mu_N^{\beta,u}(T_N > \lambda N) \to 1$.  Analogous statements hold for the annealed and/or loosened systems.  From this and \eqref{E:weightratio2} we obtain that when there is pinning in the original quenched system, i.e. for $u > u_c^q(\beta)$, we have the analog 
of \eqref{E:freeenrel}:
\[
  \beta \hf^q(\beta,u + \beta^{-1} \log M_E^f(b)) = \beta f^q(\beta,u) + b.
  \]
Therefore by continuity of $\hf^q(\beta,\cdot)$, for some $\gamma>0$, since $f^q(\beta,u) = 0$ 
for all $u < u_c(\beta,P^X)$,
\begin{align}
  u_c^q(\beta,P^X) &= \inf\{ u \in \RR: \beta \hf^q(\beta,u + \beta^{-1} \log M_E^f(b)) > b \} \notag \\
  &= \inf\{ u \in \RR: \beta \hf^q(\beta,u + \beta^{-1} \log M_E^f(b)) > 0 \} + \gamma \notag \\
  &= u_c^q(\beta,\hP^X) - \beta^{-1} \log M_E^f(b) + \gamma, \notag
  \end{align}
and the part of the free energy or contact fraction graph (as a function of $u$) with $u > u_c^q(\beta,P^X)$ in the original system is just a translate of part of the corresponding graph of the loosened system:  for all $\Delta>0$,
\begin{align} \label{E:Ccurves}
  C^q(\beta,u_c^q(\beta,P^X) + \Delta) &= 
    \hat{C}^q(\beta,u_c^q(\beta,P^X) + \beta^{-1} \log M_E^f(b) + \Delta) \\
  &= \hat{C}^q(\beta,\hu_c^q(\beta) + \gamma + \Delta) \notag \\
  &\geq \hat{C}^q(\beta,\hu_c^q(\beta) + \gamma). \notag   
  \end{align}
(Note that in contrast to the recurrent case, it is not true that the entire graph is such a translate.)  Since the quantity on the right side of \eqref{E:Ccurves}, which we now call $y$, is strictly positive, this shows that the transition in discontinuous in the original quenched system.  

The idea here is that the part of the supercritical loosened-system free-energy graph (from height $b>0$ upward) that is translated to obtain the supercritical original-system graph is bounded away from the critical point, so the contact fraction is bounded away from 0.  

To see that the quenched and annealed critical points differ, we observe that by \eqref{E:anfreeen}, \eqref{E:qfreeineq} and \eqref{E:Ccurves}, we have 
\[
  f^a(\beta,u) - f^q(\beta,u) \geq \half \beta^2 y^2 \quad \text{for all } u > u_c^q(\beta),
  \]
and the result then follows from continuity of the free energies.
\end{proof}


\begin{thebibliography}{99}

\bibitem{Al06} Alexander, K.S., \emph{The effect of disorder on polymer depinning transitions}, arXiv.org:  math.PR/0610008 (2006).

\bibitem{AS06} Alexander, K.S. and Sidoravicius, V., \emph{Pinning of polymers
and interfaces by random potentials}, Ann. Appl. Probab. \textbf{16} (2006), 636--669.

\bibitem{CGG06} Caravenna, F., Giacomin, G. and Gubinelli, M., \emph{A numerical approach to copolymers at selective interfaces}, J. Stat. Phys. \textbf{122} (2006), 799--832.

\bibitem{CCFS89}  Chayes, J. T., Chayes, L., Fisher, D. S. and Spencer, T., \emph{Correlation length bounds for disordered ferromagnets},  Commun. Math. Phys. \textbf{120} (1989), 501--523.

\bibitem{Gi06} Giacomin, G., \emph{Random Polymer Models}.  Imperial College Press, London, 2006, in press.

\bibitem{GT06a} Giacomin, G. and Toninelli, F. L., \emph{Smoothing effect of
quenched disorder on polymer depinning transitions}, Commun. Math. Phys. \textbf{266} (2006), 1-16.

\bibitem{GT06b} Giacomin, G. and Toninelli, F. L., \emph{The localized phase of disordered copolymers with adsorption}, ALEA Lat. Am. J. Probab. Math. Stat. \textbf{1} (2006), 149--180.

\bibitem{GT06c} Giacomin, G. and Toninelli, F. L., \emph{Force-induced depinning of directed polymers}, arXiv.org: cond-mat/0610663 (2006).

\bibitem{IOW04} Iliev, G., Orlandini, E. and Whittington, S. G., \emph{Adsorption and localization of random copolymers subject to a force:  The Morita approximation}, Eur. Phys. J. B \textbf{40} (2004), 63--71.

\bibitem{NN01} Naidenov, A. and Nechaev, S.,  \emph{Adsorption of a random
heteropolymer at a potential well revisited:  location of transition point and
design of sequences}, J. Phys. A: Math. Gen. \textbf{34} (2001). 5625--5634.

\end{thebibliography}
\end{document}